\documentclass[10.9pt,twoside]{amsart}
\usepackage{amsmath, amsthm, amscd, amsfonts, amssymb, graphicx, color}
\usepackage[bookmarksnumbered, colorlinks, plainpages]{hyperref}

\theoremstyle{definition}

\theoremstyle{remark}

\numberwithin{equation}{section}

\begin{document}
\setcounter{page}{1}
\begin{center}
{\bf  A GENERALIZATION OF THE WEAK AMENABILITY\\ OF SOME BANACH ALGEBRA }
\end{center}

\title[]{}
\author[]{KAZEM HAGHNEJAD AZAR  AND ABDOLHAMID RIAZI}

\address{}

\dedicatory{}



\begin{abstract} Let $A$ be a Banach algebra and $A^{**}$ be the second dual of it. We show that by some new conditions, $A$ is weakly amenable whenever $A^{**}$ is weakly amenable. We will study this problem under generalization, that is, if $(n+2)-th$ dual of $A$, $A^{(n+2)}$, is $T-S-$weakly amenable, then $A^{(n)}$ is $T-S-$weakly amenable where $T$ and $S$ are continuous linear mappings from $A^{(n)}$ into $A^{(n)}$.

\end{abstract} \maketitle

\section{\bf  Preliminaries and
Introduction }

\noindent   Let $A$ be  a Banach algebra and $A^*$,
$A^{**}$, respectively, are the first and second dual of $A$.  For $a\in A$
 and $a^\prime\in A^*$, we denote by $a^\prime a$
 and $a a^\prime$ respectively, the functionals on $A^*$ defined by $ \langle a^\prime a,b\rangle  = \langle  a^\prime,ab\rangle  =a^\prime(ab)$ and $ \langle  a a^\prime,b\rangle  = \langle  a^\prime,ba\rangle  =a^\prime(ba)$ for all $b\in A$.
   The Banach algebra $A$ is embedded in its second dual via the identification
 $ \langle  a,a^\prime\rangle  $ - $ \langle  a^\prime,a\rangle  $ for every $a\in
A$ and $a^\prime\in
A^*$. Arens [1] has shown that given any Banach algebra $A$, there exist
two algebra multiplications on the second dual  of $A$ which
extend multiplication on  $A$. In the following, we introduce both
multiplication which are given in [13]. The first (left) Arens product of $a^{\prime\prime},b^{\prime\prime}\in A^{**}$ shall be simply indicated by $a^{\prime\prime}b^{\prime\prime}$ and defined by the three steps:
 $$ \langle  a^\prime a,b\rangle  = \langle  a^\prime ,ab\rangle  ,$$
  $$ \langle  a^{\prime\prime} a^\prime,a\rangle  = \langle  a^{\prime\prime}, a^\prime a\rangle  ,$$
  $$ \langle  a^{\prime\prime}b^{\prime\prime},a^\prime\rangle  = \langle  a^{\prime\prime},b^{\prime\prime}a^\prime\rangle  .$$
 for every $a,b\in A$ and $a^\prime\in A^*$. Similarly, the second (right) Arens product of $a^{\prime\prime},b^{\prime\prime}\in A^{**}$ shall be  indicated by $a^{\prime\prime}ob^{\prime\prime}$ and defined by :
 $$ \langle  a oa^\prime ,b\rangle  = \langle  a^\prime ,ba\rangle  ,$$
  $$ \langle  a^\prime oa^{\prime\prime} ,a\rangle  = \langle  a^{\prime\prime},a oa^\prime \rangle  ,$$
  $$ \langle  a^{\prime\prime}ob^{\prime\prime},a^\prime\rangle  = \langle  b^{\prime\prime},a^\prime ob^{\prime\prime}\rangle  .$$
  for all $a,b\in A$ and $a^\prime\in A^*$.\\
  We say that $A$ is Arens regular if both multiplications are equal.
   Let $a^{\prime\prime}$ and $b^{\prime\prime}$ be elements of $A^{**}$. By $Goldstine^,s$ Theorem [6, P.424-425], there are nets $(a_{\alpha})_{\alpha}$ and $(b_{\beta})_{\beta}$ in $A$ such that $a^{\prime\prime}=weak^*-lim_{\alpha}a_{\alpha}$ ~and~  $b^{\prime\prime}=weak^*-lim_{\beta}b_{\beta}$. So it is easy to see that for all $a^\prime\in A^*$,
$$lim_{\alpha}lim_{\beta} \langle  a^\prime,a_{\alpha}b_{\beta}\rangle  = \langle  a^{\prime\prime}b^{\prime\prime},a^\prime\rangle  $$ and
$$lim_{\beta}lim_{\alpha} \langle  a^\prime,a_{\alpha}b_{\beta}\rangle  = \langle  a^{\prime\prime}ob^{\prime\prime},a^\prime\rangle  .$$
Thus $A$ is Arens regular if and only if for every $a^\prime\in A^*$, we have
$$lim_{\alpha}lim_{\beta} \langle  a^\prime,a_{\alpha}b_{\beta}\rangle  =lim_{\beta}lim_{\alpha} \langle  a^\prime,a_{\alpha}b_{\beta}\rangle  .$$
For more detail see [6, 13, 15].\\
   Let $X$ be a   Banach $A-bimodule$.
   A derivation from $A$ into $X$ is a bounded linear mapping $D:A\rightarrow X$ such that $$D(xy)=xD(y)+D(x)y~~for~~all~~x,~y\in A.$$
The space of continuous derivations from $A$ into $X$ is denoted by $Z^1(A,X)$.\\
Easy example of derivations are the inner derivations, which are given for each $x\in X$ by
$$\delta_x(a)=ax-xa~~for~~all~~a\in A.$$
The space of inner derivations from $A$ into $X$ is denoted by $N^1(A,X)$.
The Banach algebra $A$ is said to be a amenable, when for every Banach $A-bimodule$ $X$, the inner derivations are only derivations existing from $A$ into $X^*$. It is clear that $A$ is amenable if and only if $H^1(A,X^*)=Z^1(A,X^*)/ N^1(A,X^*)=\{0\}$.\\
A Banach algebra $A$ is said to be a weakly amenable, if every derivation from $A$ into $A^*$ is inner. Similarly, $A$ is weakly amenable if and only if $H^1(A,A^*)=Z^1(A,A^*)/ N^1(A,A^*)=\{0\}$.\\
 Suppose that $A$ is a Banach algebra and $X$ is a Banach $A-bimodule$. According to [5, pp.27 and 28], $X^{**}$ is a Banach $A^{**}-bimodule$, where  $A^{**}$ is equipped with the first Arens product.

\noindent Let $A^{(n)}$ and  $X^{(n)}$  be $n-th~dual$ of $A$ and $X$, respectively. By [19, page 4132-4134], if $n\geq 0$ is an even number, then  $X^{(n)}$ is a Banach $A^{(n)}-bimodule$. Then for $n\geq 2$,   we define  $X^{(n)}X^{(n-1)}$ as a subspace of $A^{(n-1)}$, that is, for all $x^{(n)}\in X^{(n)}$,  $x^{(n-1)}\in X^{(n-1)}$ and  $a^{(n-2)}\in A^{(n-2)}$ we define
$$ \langle  x^{(n)}x^{(n-1)},a^{(n-2)}\rangle  = \langle  x^{(n)},x^{(n-1)}a^{(n-2)}\rangle  .$$
If $n$ is odd number, then for $n\geq 1$,   we define  $X^{(n)}X^{(n-1)}$ as a subspace of $A^{(n)}$, that is, for all $x^{(n)}\in X^{(n)}$,  $x^{(n-1)}\in X^{(n-1)}$ and  $a^{(n-1)}\in A^{(n-1)}$ we define
$$ \langle  x^{(n)}x^{(n-1)},a^{(n-1)}\rangle  = \langle  x^{(n)},x^{(n-1)}a^{(n-1)}\rangle  .$$
If $n=0$, we take $A^{(0)}=A$ and $X^{(0)}=X$.\\
Now let $X$ be a Banach $A-bimodule$ and $D:A\rightarrow X$ be a derivation. A problem which is of interest is under what conditions  $D^{\prime\prime}$  is again a derivation. In [14, 5.9 ], this problem has been studied for the spacial case $X=A$, and they showed that  $D^{\prime\prime}$ is a derivation if and only if  $D^{\prime\prime} (A^{**})A^{**}\subseteq A^{*}$. We study this problem in the generality, that is, if  $A^{(n+2)}$ is $T-S-$weakly amenable, then it follows that $A^{(n)}$ is $T-S-$ weakly amenable where $T$ and $S$ are continuous linear mapping from $A^{(n)}$ into $A^{(n)}$ and $n\geq 0$.\\
The main results of this paper can be summarized as follows:\\
{\bf a)} Assume that $A$ is a Banach algebra  and $A^{(n+2)}$ has $T-w^*w$ property. If  ${A^{(n+2)}}$ is weakly $ T^{\prime\prime}-S^{\prime\prime}-$amenable, then ${A^{(n)}}$ is weakly $ T-S-$amenable.\\
{\bf b)} Let $X$ be  a   Banach $A-bimodule$ and let  $T,~S:A^{(n)}\rightarrow A^{(n)}$ be continuous linear mappings. Let the mapping $a^{(n+2)}\rightarrow x^{(n+2)}T^{\prime\prime}(a^{(n+2)})$ be $weak^*$-to-$weak$ continuous for all $x^{(n+2)}\in X^{(n+2)}$. Then if  $D:A^{(n)}\rightarrow X^{(n+1)}$ is a $ T-S-derivation$, it follows that $D^{\prime\prime}:A^{(n+2)}\rightarrow X^{(n+3)}$ is a $ T^{\prime\prime}-S^{\prime\prime}-derivation$.\\
{\bf c)} Let $X$  be a   Banach $A-bimodule$ and the mapping $a^{\prime\prime}\rightarrow x^{\prime\prime}a^{\prime\prime}$ be $weak^*-to- weak$ continuous for all $x^{\prime\prime}\in X^{**}$. If $D:A\rightarrow X^*$ is a derivation, then $D^{\prime\prime}(A^{**})X^{**}\subseteq A^*$.\\
{\bf d)} Let $X$ be a   Banach $A-bimodule$ and $D:A\rightarrow X^*$ be a derivation. Suppose that $D^{\prime\prime}:A^{**}\rightarrow X^{***}$ is     surjective derivation. Then  the mapping $a^{\prime\prime}\rightarrow x^{\prime\prime}a^{\prime\prime}$ is $weak^*-to- weak$ continuous for all $x^{\prime\prime}\in X^{**}$.\\
{\bf e)} Suppose that  $X$  is a   Banach $A-bimodule$ and $A$ is Arens regular. Assume that   $D:A\rightarrow X^{*}$ is a derivation  and surjective. Then $D^{\prime\prime}:A^{**}\rightarrow X^{***}$ is a derivation if and only if the mapping $a^{\prime\prime}\rightarrow x^{\prime\prime}a^{\prime\prime}$ is $weak^*-to- weak$ continuous for all $x^{\prime\prime}\in X^{**}$.\\
In every parts of this paper, $n\geq 0$ is even number.\\\\

\begin{center}
\section{ \bf Weak amenability of Banach algebras  }
\end{center}

\noindent{\it{\bf Definition 2-1.}} Let $X$ be a   Banach $A-bimodule$ and $T$, $S$ be continuous linear mappings from $A$ into itself. We say that $D:A\rightarrow X$  is $ T-S-derivation$, if
$$D(xy)=T(x)D(y)+D(x)S(y)~~for~~all~~x,~y\in A.$$
Now let $x\in A$. Then we say that the linear mapping $\delta_x:A\rightarrow A$ is inner $ T-S-derivation$, if for every $a\in A$ we have $\delta_x(a)=T(a)x-xS(a)$. \\
The Banach algebra $A$ is said to be a $T-S-$amenable, when for every Banach $A-bimodule$ $X$, every $T-S-$derivations from $A$ into $X^*$ is  inner $T-S-$derivations.\\
The definition of weakly $T-S-$ amenable is similar.\\\\
\noindent{\it{\bf Definition 2-2.}} Assume that $A$ is a Banach algebra and $T:A\rightarrow A$ is a continuous linear mapping such that the mapping  $b^{\prime\prime}\rightarrow a^{\prime\prime} T^{\prime\prime}(b^{\prime\prime})~:~A^{**}\rightarrow A^{**}~is~~~weak^*-to-weak~$ continuous where $a^{\prime\prime}\in A^{**}$. Then we say that $a^{\prime\prime}\in A^{**}$ has $T-w^*w$ property.\\
We say that $B\subseteq A^{**}$ has $T-w^*w$ property, if every $b\in B$ has $T-w^*w$ property.\\\\

 \noindent Let $A$ be a Banach algebra and $A^{**}$ has  $I-w^*w$ property whenever $I:A\rightarrow A$ is the identity mapping. Then, obviously that $A$ is Arens regular. There are some non-reflexive Banach algebras which the second dual of them  have $T-w^*w$ property.  If $A$ is Arens regular, then, in general, $A^{**}$ has not $I-w^*w$ property. In the following we give some examples from Banach algebras that the second dual of them have $T-w^*w$ property or no.
 \begin{enumerate}
\item  Let $A$ be a non-reflexive Banach space and suppose that  $ \langle  f,x\rangle  =1$ for some $f\in A^*$ and $x\in A$. We define the product on $A$ as $ab= \langle  f,a\rangle  b$ for all $a, ~b\in A$. It is clear that   $A$ is a Banach algebra with this product, then  $A^{**}$ has  $I-w^*w$ property whenever $I:A\rightarrow A$ is the identity mapping.
\item  Every reflexive Banach algebra has $T-w^*w$ property.
\item  Consider the algebra $c_0=(c_0,.)$ is the collection of all sequences of scalars that convergence to $0$, with the some vector space operations and norm as $\ell_\infty$.  Then $c_0^{**}=\ell_\infty$ has $I-w^*w$ property whenever $I:c_0\rightarrow c_0$ is the identity mapping.
\item  $L^1(G)^{**}$ and $M(G)^{**}$ have not $I-w^*w$ property whenever $G$ is locally compact group, but when $G$ is finite,  $L^1(G)^{**}$ and $M(G)^{**}$ have $I-w^*w$ property.\\\\

\end{enumerate}

\noindent{\it{\bf Theorem 2-3.}} Assume that $A$ is a Banach algebra  and $A^{(n+2)}$ has $T-w^*w$ property. If $D:A^{(n)}\rightarrow  A^{(n+1)}$ is a $ T-S-derivation$, then $D^{\prime\prime}:A^{(n+2)}\rightarrow  A^{(n+3)}$
is a $ T^{\prime\prime}-S^{\prime\prime}-derivation$.
\begin{proof}  Let $a^{(n+2)},~ b^{(n+2)}\in A^{(n+2)}$  and let $(a_\alpha^{(n)})_\alpha,~ (b_\beta^{(n)})_\beta\subseteq A^{(n)}$ such that $a_\alpha^{(n)}\stackrel{w^*} {\rightarrow}a^{(n+2)}$ and
$b_\beta^{(n)}\stackrel{w^*} {\rightarrow}b^{(n+2)}$. Due to $A^{(n+2)}$ has $T-w^*w$ property, we have
$c^{(n+2)}T(a_\alpha^{(n)})\stackrel{w} {\rightarrow}c^{(n+2)}T^{\prime\prime}(a^{(n+2)})$ for all $c^{(n+2)}\in A^{(n+2)}$. Using the $weak^*-to-weak^*$ continuity of $D^{\prime\prime}$, we obtain
 $$lim_\alpha lim_\beta  \langle  T(a_\alpha^{(n)})D(b_\beta^{(n)}),c^{(n+2)}\rangle   = lim_\alpha lim_\beta \langle  D(b_\beta^{(n)}),c^{(n+2)}T(a_\alpha^{(n)})\rangle   $$$$=lim_\alpha  \langle  D^{\prime\prime}(b^{(n+2)}),c^{(n+2)}T(a_\alpha^{(n)})\rangle   =
  \langle  D^{\prime\prime}(b^{(n+2)}),c^{(n+2)}T^{\prime\prime}(a^{(n+2)})\rangle   $$
 $$=  \langle  T^{\prime\prime}(a^{(n+2)})D^{\prime\prime}(b^{(n+2)}),c^{(n+2)}\rangle  .$$
 Moreover, it is also clear that for every $c^{(n+2)}\in A^{(n+2)}$, we have
 $$lim_\alpha lim_\beta  \langle  D(a_\alpha^{(n)})S(b_\beta^{(n)}),c^{(n+2)}\rangle   = \langle  D^{\prime\prime}(a^{(n+2)})S^{\prime\prime}(b^{(n+2)}),c^{(n+2)}\rangle  .$$
 Notice that in latest equalities, we didn't need $S-w^*w$ property for $A^{(n+2)}$.
 In the following, we take limit on the $weak^*$ topologies. Thus we have
 $$D^{\prime\prime}( a^{(n+2)}b^{(n+2)})=lim_\alpha lim_\beta D(a_\alpha^{(n)}b_\beta^{(n)})=
lim_\alpha lim_\beta T(a_\alpha^{(n)})D(b_\beta^{(n)})+$$
$$lim_\alpha lim_\beta D(a_\alpha^{(n)})S(b_\beta^{(n)})=
 T^{\prime\prime}(a^{(n+2)})D^{\prime\prime}(b^{(n+2)})+D^{\prime\prime}( a^{(n+2)})S^{\prime\prime}(b^{(n+2)}).$$
\end{proof}

\noindent{\it{\bf Theorem 2-4.}}  Assume that $A$ is a Banach algebra  and $A^{(n+2)}$ has $T-w^*w$ property. If  ${A^{(n+2)}}$ is weakly $ T^{\prime\prime}-S^{\prime\prime}-$amenable, then ${A^{(n)}}$ is weakly $ T-S-$amenable.
\begin{proof} Let $D:A^{(n)}\rightarrow  A^{(n+1)}$ is a $ T-S-derivation$, then by Theorem 2-3, $D^{\prime\prime}:A^{(n+2)}\rightarrow  A^{(n+3)}$
is a $ T^{\prime\prime}-S^{\prime\prime}-derivation$. Since $A^{(n+2)}$ is weakly $ T^{\prime\prime}-S^{\prime\prime}-$ amenable, $D^{\prime\prime}:A^{(n+2)}\rightarrow  A^{(n+3)}$
is an inner $ T^{\prime\prime}-S^{\prime\prime}-~~derivation$. It follows that for every $a^{(n+2)}\in A^{(n+2)}$, we have
$$D^{\prime\prime}(a^{(n+2)})=T^{\prime\prime}(a^{(n+2)})a^{(n+3)}-a^{(n+3)}S^{\prime\prime}(a^{(n+2)}).$$
for some $a^{(n+3)}\in A^{(n+3)}$.
Take $a^{(n+1)}=a^{(n+3)}\mid_{A^{(n+1)}}$. Then for every $a^{(n)}\in A^{(n)}$, we have
$$D(a^{(n)})=T(a^{(n)})a^{(n+1)}-a^{(n+1)}S(a^{(n)}).$$
It follows that $D$ is inner $T-S-$derivation, and so proof is hold.

\end{proof}

\noindent{\it{\bf Corollary 2-5.}} Let $A$ be a Banach algebra  and  $I:A\rightarrow A$ be identity mapping. If $A^{**}$ has $I-w^*w$ property and $A^{**}$ is weakly amenable, then $A$ is weakly amenable.\\

\noindent{\it{\bf Corollary 2-6.}} Let $A$ be a Banach algebra. If $A^{***}A^{**}\subseteq A^{*}$ and $A^{**}$ is weakly amenable, then $A$ is weakly amenable.
\begin{proof} We show that $A^{**}$ has $I-w^*w$ property where $I:A\rightarrow A$ is identity mapping. Suppose that $a^{\prime\prime},~b^{\prime\prime}\in A^{**}$ and $b_{\alpha}^{\prime\prime}\stackrel{w^*} {\rightarrow}b^{\prime\prime}$. Let $c^{\prime\prime\prime}\in A^{***}$. Since $c^{\prime\prime\prime}a^{\prime\prime}\in A^*$, we have
$$ \langle  c^{\prime\prime\prime},a^{\prime\prime}b_{\alpha}^{\prime\prime}\rangle
= \langle  c^{\prime\prime\prime}a^{\prime\prime},b_{\alpha}^{\prime\prime}\rangle  =
 \langle  b_{\alpha}^{\prime\prime},c^{\prime\prime\prime}a^{\prime\prime}\rangle  \rightarrow
 \langle  b^{\prime\prime},c^{\prime\prime\prime}a^{\prime\prime}\rangle  =
 \langle  c^{\prime\prime\prime},a^{\prime\prime}b^{\prime\prime}\rangle  .$$
We conclude that $a^{\prime\prime}b_{\alpha}^{\prime\prime}\stackrel{w} {\rightarrow}a^{\prime\prime}b^{\prime\prime}$. So $A^{**}$ has $I-w^*w$ property. By using Corollary 2-5, $A$ is weakly amenable.\\\end{proof}

\noindent{\it{\bf Example 2-7.}} $c_0$ is weakly amenable.
\begin{proof} Since $\ell^\infty=c_0^{**}$ is weakly amenable and $\ell^\infty$ has $I-w^*w$ property by Corollary 2-5, proof is hold.\\
\end{proof}

\noindent{\it{\bf Theorem 2-8.}} Suppose that  $A$ is a Banach algebra and  $B$ is a closed subalgebra of $A^{(n+2)}$ that is consisting of $A^{(n)}$ where $n\in \mathbb{N}\cup \{0\}$. If $B$ has $T-w^*w$ property and  is weakly $ T^{\prime\prime}-S^{\prime\prime}-$amenable, then $A^{(n)}$ is weakly $ T-S-$amenable.

\begin{proof} Suppose that $D:A^{(n)}\rightarrow A^{(n+1)}$ is a $ T-S-derivation$ and $p:A^{(n+3)}\rightarrow B^\prime$ is the restriction map, defined by $P(a^{(n+3)})=a^{(n+3)}\mid_{B^\prime}$ for every $a^{(n+3)}\in A^{(n+3)}$. Since $B$  has $T-w^*w$ property, $\bar{D}=PoD^{\prime\prime}\mid_B:B\rightarrow B^\prime$ is a $ T^{\prime\prime}-S^{\prime\prime}-derivation$. Since $B$ is weakly $ T^{\prime\prime}-S^{\prime\prime}-$amenable, there is $b^\prime\in B^\prime$ such that $\bar{D}=\delta_{b^\prime}$. We take $a^{(n+1)}=b^\prime\mid_{A^{(n+1)}}$, then $D=\bar{D}$ on $A^{(n+1)}$. Consequently, we have $D=\delta_{a^{(n+1)}}$.

\end{proof}

\noindent{\it{\bf Theorem 2-9.}} Let $X$ be a   Banach $A-bimodule$ and let  $T,~S:A^{(n)}\rightarrow A^{(n)}$ be continuous linear mappings. Let the mapping $a^{(n+2)}\rightarrow x^{(n+2)}T^{\prime\prime}(a^{(n+2)})$ be $weak^*$-to-$weak$ continuous for all $x^{(n+2)}\in X^{(n+2)}$. If  $D:A^{(n)}\rightarrow X^{(n+1)}$ is a $ T-S-derivation$, then $D^{\prime\prime}:A^{(n+2)}\rightarrow X^{(n+3)}$ is a $ T^{\prime\prime}-S^{\prime\prime}-derivation$.

\begin{proof}
Let $a^{(n+2)},~ b^{(n+2)}\in A^{(n+2)}$  and let $(a_\alpha^{(n)})_\alpha,~ (b_\beta^{(n)})_\beta\subseteq A^{(n)}$ such that $a_\alpha^{(n)}\stackrel{w^*} {\rightarrow}a^{(n+2)}$ and
$b_\beta^{(n)}\stackrel{w^*} {\rightarrow}b^{(n+2)}$. Then for all $x^{(n+2)}\in X^{(n+2)}$, we have
 $x^{(n+2)}T(a_\alpha^{(n)})\stackrel{w} {\rightarrow}x^{(n+2)}a^{(n+2)}$. Consequently, we have
$$lim_\alpha lim_\beta \langle  T(a_\alpha^{(n)})D(b_\beta^{(n)}),x^{(n+2)}\rangle  =lim_\alpha lim_\beta \langle  D(b_\beta^{(n)}),x^{(n+2)}T(a_\alpha^{(n)})\rangle  $$$$=lim_\alpha \langle  D^{\prime\prime}(b^{(n+2)}),x^{(n+2)}T(a_\alpha^{(n)})\rangle  = \langle  D^{\prime\prime}(b^{(n+2)}),x^{(n+2)}T(a^{(n+2)})\rangle  $$$$=
 \langle  T(a^{(n+2)})D^{\prime\prime}(b^{(n+2)}),x^{(n+2)}\rangle  .$$
For every $x^{(n+2)}\in X^{(n+2)}$, we have also the following equalities
$$lim_\alpha lim_\beta \langle  D(a_\alpha^{(n)})S(b_\beta^{(n)}),x^{(n+2)}\rangle  =lim_\alpha lim_\beta \langle  D(a_\alpha^{(n)}),S(b_\beta^{(n)})x^{(n+2)}\rangle  $$$$=lim_\alpha  \langle  D(a_\alpha^{(n)}),S(b^{(n+2)})x^{(n+2)}\rangle  =
 \langle  D^{\prime\prime}(a^{(n+2)}),S(b^{(n+2)})x^{(n+2)}\rangle  $$$$= \langle  D^{\prime\prime}(a^{(n+2)})S(b^{(n+2)}),x^{(n+2)}\rangle  .$$
 In the following, we take limit on the $weak^*$ topologies. Using the $weak^*-to-weak^*$ continuity of $D^{\prime\prime}$, we obtain
$$D^{\prime\prime}(a^{(n+2)}b^{(n+2)})=lim_\alpha lim_\beta D(a_\alpha^{(n)}b_\beta^{(n)})=
lim_\alpha lim_\beta T(a_\alpha^{(n)})D(b_\beta^{(n)})+$$$$lim_\alpha lim_\beta D(a_\alpha^{(n)})S(b_\beta^{(n)})=
T^{\prime\prime}(a^{(n+2)})D^{\prime\prime}(b^{(n+2)})+D^{\prime\prime}(a^{(n+2)})S^{\prime\prime}(b^{(n+2)}).$$
Thus $D^{\prime\prime}:A^{(n+2)}\rightarrow X^{(n+3)}$ is a $ T^{\prime\prime}-S^{\prime\prime}-derivation$.

\end{proof}

\noindent{\it{\bf Corollary 2-10.}} Let $X$ be a   Banach $A-bimodule$ and the mapping $a^{\prime\prime}\rightarrow x^{\prime\prime}a^{\prime\prime}$ be $weak^*-to- weak$ continuous for all $x^{\prime\prime}\in X^{**}$. Then, if
$H^1(A^{**},X^{***})={0}$, it follows that $H^1(A,X^*)={0}$.\\\\

\noindent{\it{\bf Corollary 2-11.}} Let $X$  be a   Banach $A-bimodule$ and the mapping $a^{\prime\prime}\rightarrow x^{\prime\prime}a^{\prime\prime}$ be $weak^*-to- weak$ continuous for all $x^{\prime\prime}\in X^{**}$. If $D:A\rightarrow X^*$ is a derivation, then $D^{\prime\prime}(A^{**})X^{**}\subseteq A^*$.

\begin{proof} By using Theorem 2-9 and [14, Corollary 4-3], proof is hold.

\end{proof}

\noindent{\it{\bf Theorem 2-12.}} Let $X$ be a   Banach $A-bimodule$ and $D:A\rightarrow X^*$ be a  surjective derivation. Suppose that $D^{\prime\prime}:A^{**}\rightarrow X^{***}$ is  also a   derivation. Then  the mapping $a^{\prime\prime}\rightarrow x^{\prime\prime}a^{\prime\prime}$ is $weak^*-to- weak$ continuous for all $x^{\prime\prime}\in X^{**}$.

\begin{proof} Let $a^{\prime\prime}\in A^{**}$ such that $a_{\alpha}^{\prime\prime}\stackrel{w^*} {\rightarrow}a^{\prime\prime}$. We show that $x^{\prime\prime}a_{\alpha}^{\prime\prime}\stackrel{w} {\rightarrow}x^{\prime\prime}a^{\prime\prime}$  for all $x^{\prime\prime}\in X^{**}$. Suppose that $x^{\prime\prime\prime}\in X^{***}$. Since $D^{\prime\prime}(A^{**})=X^{***}$, by using [14, Corollary 4-3], we conclude that  $X^{***}X^{**}=D^{\prime\prime}(A^{**})X^{**}\subseteq A^*$. Then  $x^{\prime\prime\prime}x^{\prime\prime}\in A^*$, and so we have the following equality
$$ \langle  x^{\prime\prime\prime},x^{\prime\prime}a_\alpha^{\prime\prime}\rangle  = \langle  x^{\prime\prime\prime}x^{\prime\prime},a_\alpha ^{\prime\prime}\rangle  = \langle  a_\alpha ^{\prime\prime},x^{\prime\prime\prime}x^{\prime\prime}\rangle  \rightarrow
 \langle  a^{\prime\prime},x^{\prime\prime\prime}x^{\prime\prime}\rangle  =
 \langle  x^{\prime\prime\prime},x^{\prime\prime}a^{\prime\prime}\rangle  .$$

\end{proof}

\noindent{\it{\bf Corollary 2-13.}} Suppose that  $X$  is a   Banach $A-bimodule$ and $A$ is Arens regular. Assume that   $D:A\rightarrow X^{*}$ is a surjective derivation. Then $D^{\prime\prime}:A^{**}\rightarrow X^{***}$ is a derivation if and only if the mapping $a^{\prime\prime}\rightarrow x^{\prime\prime}a^{\prime\prime}$ from $A^{**}$ into $X^{**}$ is $weak^*-to- weak$ continuous for all $x^{\prime\prime}\in X^{**}$.

\begin{proof} By using Corollary 2-11, Theorem 2-12 and [14, Corollary 4-3], proof is hold.

\end{proof}

 In the proceeding Corollary, if we  omit the Arens regularity of $A$, then we  have also the following conclusion.\\
Assume that   $D:A\rightarrow X^{*}$ is a surjective derivation. Then, $D^{\prime\prime}(A^{**})X^{**}\subseteq A^*$ if and only if the mapping $a^{\prime\prime}\rightarrow x^{\prime\prime}a^{\prime\prime}$ is $weak^*-to- weak$ continuous for all $x^{\prime\prime}\in X^{**}$.\\

\noindent \noindent{\it{\bf Corollary 2-14.}}  Let $A$ be a Banach algebra. Then we have the following results:
\begin{enumerate}
\item   Assume that  $A$ is Arens regular and $D:A\rightarrow A^*$ is a surjective derivation. Then $D^{\prime\prime}:A^{**}\rightarrow A^{***}$ is a derivation if and only if $A$ has $I-w^*w$ property whenever $I:A\rightarrow A$ is the identity mapping.
\item  Assume that   $D:A\rightarrow A^{*}$ is a surjective derivation. Then,  $A$ has $I-w^*w$ property if and only if $D^{\prime\prime}(A^{**})A^{**}\subseteq A^*$. So it is clear that if $D:A\rightarrow A^{*}$ is a surjective derivation and $D^{\prime\prime}(A^{**})A^{**}\subseteq A^*$, then $A$ is Arens regular.\\
\end{enumerate}

\noindent{\it{\bf Problem.}} Let $S$ be a semigroup. Dose $C(S)^{**}$, $L^1(S)^{**}$ and $M(S)^{**}$ have $I-w^*w$ property? whenever $I$ is the identity mapping.\\\\

\bibliographystyle{amsplain}

\begin{thebibliography}{99}
\bibitem{1} R. E.  Arens, {\it The adjoint of a bilinear operation}, Proc. Amer. Math. Soc. {\bf 2} (1951), 839-848.
\bibitem{2} N. Arikan, {\it A simple condition ensuring the Arens
regularity of bilinear mappings}, Proc. Amer. Math. Soc. {\bf 84}
(4) (1982), 525-532.
\bibitem{3} J. Baker, A.T. Lau, J.S. Pym {\it Module homomorphism and topological centers associated with
 weakly sequentially compact Banach algebras}, Journal of Functional Analysis. {\bf 158} (1998), 186-208.
\bibitem{4} F. F. Bonsall, J. Duncan, {\it Complete normed algebras}, Springer-Verlag, Berlin 1973.
\bibitem{5} H. G. Dales, A. Rodrigues-Palacios, M.V. Velasco, {\it The second transpose of a derivation}, J. London. Math. Soc. {\bf2} 64 (2001) 707-721.
  \bibitem{6} H. G. Dales, {\it Banach algebra and automatic continuity}, Oxford 2000.

\bibitem{7} N. Dunford, J. T. Schwartz, {\it Linear operators.I},
Wiley, New york 1958.
\bibitem{8} M. Eshaghi Gordji, M. Filali, {\it Arens regularity of module actions},  Studia Math. {\bf 181} 3 (2007), 237-254.
\bibitem{8} M. Eshaghi Gordji, M. Filali, {\it Weak amenability of the second dual of a Banach algebra}, Studia Math. {\bf182} 3 (2007), 205-213.







\bibitem{9} E. Hewitt, K. A. Ross,  {\it Abstract harmonic analysis}, Springer, Berlin, Vol I 1963.
\bibitem{10} E. Hewitt, K.  A. Ross,  {\it Abstract harmonic analysis}, Springer, Berlin, Vol II 1970.
\bibitem{11}  A. T. Lau, V. Losert, {\it On the second Conjugate Algebra of locally
compact groups}, J. London Math. Soc.  {\bf 37} (2)(1988),
464-480.
\bibitem{12} A. T. Lau, A. Ulger, {\it Topological center of certain dual
algebras}, Trans. Amer.  Math. Soc. {\bf 348} (1996), 1191-1212.

\bibitem{13} S. Mohamadzadih, H. R. E. Vishki, {\it Arens regularity of module actions and the second adjoint of a derivation}, Bulletin of the Australian Mathematical Society {\bf77} (2008), 465-476.



  \bibitem{17}  J. S. Pym, {\it The convolution of functionals on spaces of bounded functions},
 Proc. London Math Soc.  {\bf 15} (1965), 84-104.

\bibitem{19} A. Ulger, {\it Some stability properties of Arens regular bilinear operators}, Proc. Amer. Math. Soc. (1991) {\bf 34}, 443-454.
\bibitem{20}  A. Ulger, {\it Arens regularity of weakly sequentialy compact Banach algebras},
Proc. Amer. Math. Soc. {\bf 127} (11) (1999), 3221-3227.
\bibitem{21} P. K. Wong, {\it The second conjugate algebras of
Banach algebras}, J. Math. Sci. {\bf 17} (1) (1994), 15-18.

\bibitem{22} Y. Zhing,  {\it Weak amenability of module extentions of Banach algebras}, Trans. Amer. Math. Soc. {\bf 354} (10) (2002), 4131-4151.
\bibitem{22} Y. Zhing,  {\it Weak amenability of a class of Banach algebra}, Cand. Math. Bull. {\bf 44} (4) (2001) 504-508.\\
\end{thebibliography}

\it{Department of Mathematics, Amirkabir University of Technology, Tehran, Iran\\
{\it Email address:} haghnejad@aut.ac.ir\\\\
Department of Mathematics, Amirkabir University of Technology, Tehran, Iran\\
{\it Email address:} riazi@aut.ac.ir}
\end{document}